# Elliptic Curves and Hyperdeterminants in Quantum Gravity
Philip Gibbs

**Abstract**

Hyperdeterminants are generalizations of determinants from matrices to multi-dimensional hypermatrices. They were discovered in the 19[th] century by Arthur Cayley but were largely ignored over a period of 100 years before once again being recognised as important in algebraic geometry, physics and number theory. It is shown that a cubic elliptic curve whose Mordell-Weil group contains a $Z_2$ x $Z_2$ x $Z$ subgroup can be transformed into the degree four hyperdeterminant on a 2x2x2 hypermatrix comprising its variables and coefficients. Furthermore, a multilinear problem defined on a 2x2x2x2 hypermatrix of coefficients can be reduced to a quartic elliptic curve whose J-invariant is expressed in terms of the hypermatrix and related invariants including the degree 24 hyperdeterminant. These connections between elliptic curves and hyperdeterminants may have applications in other areas including physics.

## Motivation

Thirty years ago I became interested in an old problem of Diophantus who had asked for sets of rational numbers such that the product of any two is one less than a square. I found that all known integer quadruples of this type satisfy a mysterious polynomial equation and I wondered where its unusual factorization properties came from. It took about twenty years before I noticed the hidden $SL(2)^3$ symmetry of the problem and learnt that the polynomial is a hyperdeterminant [1]. Thus the Diophantine problem makes interesting connections between algebraic invariants and elliptic curves. Curiously, another (possibly related) relationship between hyperdeterminants and number theory was discovered at about the same time by Manjul Bhargava [2].

In the intervening years my main interests had turned to theoretical physics, so it was good to learn that hyperdeterminants have also become important in quantum information [3] and string theory [4]. It is possible that the connections with number theory may also be relevant to physics, especially string theory where hyperdeterminants, elliptic curves and the J-function make an appearance in connection with the entropy of black holes. The objective of this paper is to suggest directions for this idea.

Since this article may be of interest to people with limited knowledge of number theory and algebraic invariants I will begin by summarising some of the concepts to be used.

## Elliptic Curves

An elliptic curve is an equation in two unknowns $x$ and $y$ of the form

$$y^2 = x^3 + \alpha x + \beta$$

This equation is defined over the complex numbers, but in Diophantine analysis, solutions over other fields and rings are of interest, especially the rationals, Integers and finite fields. The two coefficients $\alpha$ and $\beta$ are assumed to be integers.

There is a trick for constructing a new rational solution to an elliptic curve given two known ones $A = (x_1, y_1)$ and $B = (x_2, y_2)$. Picture the curve drawn on the $(x, y)$ plane and draw a straight line through the points $A$ and $B$. The line has an equation $y = mx + c$ which can be substituted into the elliptic curve to give a cubic whose roots are the $x$-coordinates of the points of intersection of the straight line and the curve. Since we already know two rational roots we can deduce that the third root of the cubic is also rational because the product of the roots is $c^2 - \beta$. This provides a third rational solution $C = (x_3, y_3)$ except in exceptional cases where it is equal to one of the known solutions. Another more trivial way to obtain a new solution $D$ is to reverse the sign of $y$ in the point $A$. These operations can be used recursively to construct more solutions.

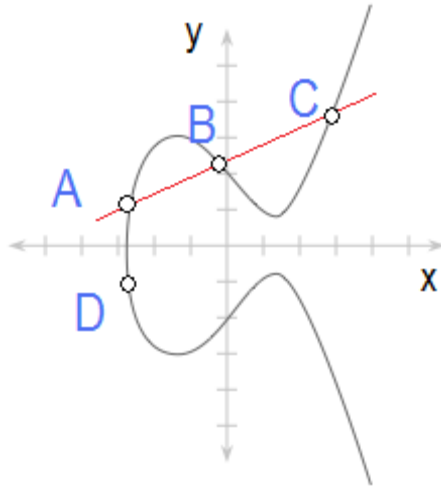

It turns out that these operations actually form an abelian group of solutions (known as the Mordell-Weil group) where the three points on the straight line satisfy the group relations $ABC = 1$ and changing the sign of $y$ is inversion, $AD = 1$. The group takes the form $T \times Z^r$ where $T$ is a finite group known as the torsion and $r$ is called the rank of the elliptic curve. $Z$ is the integers.

The details of these results can be found in any standard text book on elliptic curves.

## Cayley's Hyperdeterminant

In the nineteenth century when matrices and their determinants were newly discovered, some mathematicians wondered if the concept could be generalised to multi-dimensional hypermatrices. To do this you need a definition of the determinant that will work for more dimensions. One such definition is that the determinant of a matrix M is a discriminant that tells you when the bilinear form

$$M(u, v) = u^T M v$$

has a singular point (other than at the origin) where all its derivatives with respect to $u$ and $v$ are zero. This is true iff $u$ and $v$ are left and right eignevectors with zero eignevalue which in turn is possible iff $det M = 0$.

For a hypermatrix $A$ with three dimensions we can define a trilinear form $A(t, u, v)$. The hyperdeterminant is a multi-variable polynomial in the components of $A$ which is zero iff this trilinear form has a singular point where all derivatives are zero. The definition can be extended to hypermatrices with more dimensions.

In general it is hard to derive the form of a hyperdeterminant for a hypermatrix of given dimensions, but for simple cases such as the 2x2x2 hypermatrices the hyperdeterminant is just the discriminant of the quadratic formed by the determinant of the hypermatrix acting on one vector. The result for the hypermatrix

$$A = \left(\begin{pmatrix} a & b \\ c & d \end{pmatrix} \begin{pmatrix} e & f \\ g & h \end{pmatrix}\right)$$

Is given by the quartic polynomial

$$\det(A) = (ah + de - cf - bg)^2 - 4(ad - bc)(eh - fg)$$

## Elliptic Curves and Cayley's Hyperdeterminant

An Elliptic Curve over the rational numbers is often written as an equation between a square and a general cubic which can be reduced to the standard form above.

$$y^2 = ax^3 + bx^2 + cx + d$$

As previously stated, the solutions have a structure that forms an abelian group known as the Mordell-Weil group. In this section we will show that if the Mordell-Weil group contains a subgroup isomorphic to $Z_2$ x $Z_2$ x $Z$ then this elliptic curve can be transformed to the form of Cayley's Hyperdeterminant. In order to do this we search for hidden symmetries in the elliptic curve. As written above there is very little apparent symmetry. We can change the sign of $y$

$$y \to -y$$

or we can rescale $x$

$$x \to sx$$

$$a \to as^{-3}, \; b \to bs^{-2}, \; c \to cs^{-1}$$

To find more symmetry we need to recast the equation in a less familiar form

$$u^2v^2 = 2euv - 2fv - 2gu + h$$

Where $u$ and $v$ are regarded as the unknowns. This form for an elliptic curve is preferred because it has a $Z_2$ symmetry under an interchange,

$$u \leftrightarrow v, \; f \leftrightarrow g$$

The standard form can be retrieved by viewing this as a quadratic in $u$ and completing the square

$$(uv^2 - ev + g)^2 = -2fv^3 + (h + e^2)v^2 - 2egv + g^2$$

And substituting

$$y = uv^2 - ev + g, \; x = v$$

$$a = -2f, \quad b = h + e^2, \quad c = -2eg, \quad d = g^2$$

As a transformation over the rationals this can be inverted provided $d$ is a rational squared which is equivalent to saying that the elliptic curve has a solution at $x = 0$. We can transform any elliptic curve with at least one rational solution to this form by offsetting $x$. Since we are assuming that the group of solutions contains $Z_2 \times Z_2 \times Z$ we can take this solution to be a generator of the $Z$ subgroup. We then know that there are two more distinct solutions that generate the separate $Z_2$ subgroups. These are equal to their own inverse, but inversion in the group of solutions corresponds to reversing the sign of $y$. This means that these two solutions exist on $y = 0$, so the cubic has two rational roots. If it has two roots then it factorises fully and the third root corresponds to the product of the other two in the group of solutions.

This means we can reform the equation by factorizing the cubic

$$y^2 = 4(l - kx)(n - mx)(q - px)$$

The elliptic equation now also has an extra $S_3$ symmetry under permutations of the coefficients $k, m, p$ and $l, n, q$.

We know that $d = 4lnq$ is square. A general solution for this which preserves the $S_3$ symmetry is the substitution

$$l = rs, \quad n = ts, \quad q = rt$$

$$a = -4kmp, \quad b = 4kmrt + 4kpts + 4mprs, \quad c = -4rts(kt + mr + ps), \quad d = (2rst)^2$$

Now we can invert the transformation to find the coefficients

$$g = 2rst$$

$$e = \frac{-c}{2g} = kt + mr + ps$$

$$h = b - e^2 = 2kmrt + 2kptd + 2mprs - k^2t^2 - m^2r^2 - p^2s^2$$

$$f = \frac{-a}{2} = 2kmp$$

So the elliptic curve now takes the form

$$u^2v^2 + k^2t^2 + m^2r^2 + p^2s^2 - 2ktuv - 2mruv - 2psuv - 2kmrt - 2kptd - 2mprs + 4kmpv + 4rstu = 0$$

We originally constructed this to have a $Z_2$ symmetry and an $S_3$ symmetry but the equation now has a bigger symmetry which can be seen by placing the eight variables at the corner of a cube. The equation is symmetric under permutations of the variables generated by any rotations of the cube. So the symmetry group is $A_4$. In fact the symmetry is bigger because the equation takes the form of Cayley's hyperdeterminant which has an additional $SL(2) \times SL(2) \times SL(2)$ symmetry.

## Schläfli's Hyperdeterminant

As well as Cayley's 2x2x2 hyperdeterminant, the 2x2x2x2 hyperdeterminant (known as Schläfli's hyperdeterminant) is also linked to elliptic curves in an interesting way

Schläfli's hyperdeterminant for a 2x2x2x2 hypermatrix with components $a_{ijkl}$ is constructed by first reducing it it a 2x2x2 hypermatrix by contracting with one vector $x_l$

$$b_{ijk} = a_{ijkl}x_l$$

Cayley's hyperdeterminant for this hypermatrix then takes the form of a homogenious quartic in the two variables $x = x_0$ and $y = x_1$

$$Q(x,y) = Ax^4 + Bx^3y + Cx^2y^2 + Dxy^3 + Ey^4$$

Where $A, B, C, D, E$ are degree 4 polynomials in the components of the original hypermatrix. Schlafli's hyperdeterminant is the discriminant of this quartic which is a degree six expression in the coefficients given by

$$\Delta(Q) = S^3 - 27T^2$$

$$S = AE - 4BD + 3C^2$$

$$T = ACE + 2BCD - AD^2 - C^3 - EB^2$$

The hyperdeterminant is therefore a degree 24 polynomial in the components of the hypermatrix. The quantities $S$ and $T$ are also of interest because they are also invariants under $SL(2)^4$ transformations of the hypermatrix. $S$ is a degree 8 invariant and $T$ is degree 12.

## Elliptic curves from Schläfli's Hyperdetermiant

To see the link between Schläfli's Hyperdeterminant and elliptic curves consider the problem of finding solutions to a pair of trilinear equations

$$a_{ijkl}z_jy_kx_l = 0$$

The coefficients $a_{ijkl}$ are given integers in a 2x2x2x2 hypermatrix and the six vectors are unknowns with integer coefficients. To solve this problem we remove one vector at a time. Define a matrix $M$ with components

$$M_{ij} = a_{ijkl}y_kx_l$$

Then we require that $M_{ij}z_j = 0$ which is solvable iff $\det(M) = 0$. $M$ can now be considered as a homogeneous quadratic in the two components of $y_k$ whose coefficients depend on the components of the 2x2x2 hypermatrix with components $b_{ijk} = a_{ijkl}x_l$. This has rational solutions iff the discriminant of the quadratic is a perfect square, but this discriminant is Cayley's hyperdterminant so the requirement is that the quartic polynomial $Q(x,y)$ is a perfect square.

Although we defined an elliptic equation to be a cubic set to a square, we can also regard the quartic problem as elliptic. In fact it can be reduced to the cubic form provided it has at least one rational solution. The problem of solving the trilinear equations over the rationals is therefore equivalent to finding solutions to elliptic curves.

## The J-invariant and J-function

From the theory of elliptic curves we know that there is an important function known as the J-invariant which is independent of the rational transformations that can be applied to the elliptic curve. In terms of the invariants we have defined, the J-invariant is given by

$$J = \frac{S^3}{\Delta}$$

This expression therefore links Schläfli's hyperdeterminant and the octic invariant to an elliptic curve defined on the hypermatrix in a natural way.

The J-invariant is also tied to the J-function which is a modular form $J(\tau)$ given by the J-invariant as a function of a parameter of the solution space. The hyperdeterminant $\Delta$ is then related in a similar way to the Dedekind eta function which is also a modular form

$$\Delta(\tau) = (\eta(\tau))^{24}$$

Since the hyperdeterminant is of degree 24 this links the degree to the appearance of the number 24 in elliptic curves which is known to be related to the number of dimensions in bosonic string theory and the Leech Lattice [5]. This makes a striking connection between hyperdeterminants and other areas of mathematics and physics.

## Hyperdeterminants, Exceptional Groups and Physics

Hyperdeterminants are of interest in physics because of their relationship to the entanglement of qubits [3] and also to the entropy of STU black holes in string theory. The entropy in the STU model of 4D black holes was identified as the square root of Cayley's Hyperdeterminant by Duff [4]. After dimensional reduction in the time dimension a further $SL(2)$ symmetry is added so that the invariants of the 2x2x2x2 hypermatrix become relevant [6].

Hyperdeterminants are also of interest because of their connection with invariants of Lie groups. The cartan quartic invariant of $E_7{}^7$ on the 56 dimensional representation is a combinations of 7 2x2x2 hyperdeterminants and other invariant terms. This arises from M-theory as the squared entropy of black holes with U-duality in 4-spacetime dimensions [7].

For black holes in D=3 from dimensionally reduced supergravity, the corresponding duality group is $E_8{}^8$ [8]. The Casimir invariants of this group on the fundamental representation of dimension 248 are of degree 2, 8, 12, 14, 18, 20, 24, 30. These are independent generators of the ring of invariants [9]. It has been suggested that the associated invariant tensors for the higher orders should be constructible from just the quadratic and octic invariants [10]. When the group $SL(2)^8$ is embedded in $E_8{}^8$, the representation splits into 14 2x2x2x2 hypermatrices and 8 $SL(2)$ reps of dimension 3 (248 = 14x16+8x3) [11]. It follows that the Casimir invariants will reduce to invariants of the 2x2x2x2 hypermatrix. The primitive invariants of $SL(2)^4$ are of degree 2,4,4,6, [12] but the E8 invariants must reduce to combinations of these that are symmetric under permutations of the axis. This includes the quadratic invariant and the degree 8 and 12 invariants $S$ and $T$ given above. It is therefore conjectured that the degree 8 and 12 invariants of $E_8{}^8$ could reduce to $S$ and $T$ while the degree 24 invariant may reduce to Schläfli's hyperdeterminant. However the relationship $\Delta = S^3 - 27T^2$ cannot extend to the corresponding $E_8{}^8$ invariants because they should be independent.

In the different context of BTZ black holes in a theory of quantum gravity for 3D spacetime, the entropy of black holes may be described by the J-function [13]. Although this is a quite different model from STU back holes, the possible appearance of both the J-function and Schläfli's hyperdeterminant to describe entropy of black holes in 3D spacetime suggests a tantalising connection.

## Acknowledgments

Grateful thanks go to Leron Borsten, Martin Cederwall, Lawrence Crowell**,** Duminda Dahanayake, Mike Duff, Jakob Palmkvist, Michael Rios, Marni Sheppeard and others for discussions, especially concerning the last section of this paper.